\documentclass[11pt]{article}

\usepackage{amsthm}
\usepackage{amsmath}
\usepackage{amssymb}
\usepackage[english]{babel}
\usepackage[center]{caption2}
\usepackage{amsfonts,amssymb,amsmath,latexsym,amsthm}
\usepackage{multirow}
\usepackage[usenames,dvipsnames]{pstricks}
\usepackage{epsfig}

\newtheorem{dfn}{Definition}
\newtheorem{thm}[dfn]{Theorem}

\newtheorem{corollary}[dfn]{Corollary}
\newtheorem{conjecture}[dfn]{Conjecture}

\newtheorem{prob}[dfn]{Problem}

\def\ex{\mathrm{ex}}

\let\svthefootnote\thefootnote
\newcommand\blankfootnote[1]{%
	\let\thefootnote\relax\footnotetext{#1}%
	\let\thefootnote\svthefootnote%
}

\begin{document}
\title{On extremal numbers of the triangle plus the four-cycle} %Improving the maximum number of edges in graphs of girth at least five
\author{
Jie Ma\thanks{School of Mathematical Sciences, University of Science and Technology of China, Hefei, Anhui 230026, China.
Partially supported by the National Key R and D Program of China 2020YFA0713100,
National Natural Science Foundation of China grant 12125106, and Anhui Initiative in Quantum Information Technologies grant AHY150200.
Email: jiema@ustc.edu.cn.}
~~~~~~~
Tianchi Yang\thanks{Department of Mathematics, National University of Singapore, 119076, Singapore.
Partially supported by a start-up grant at NUS and an MOE Academic Research Fund (AcRF) Tier 1 grant. Email: tcyang@nus.edu.sg.}
}
\date{}

%\blankfootnote{}

\maketitle

\begin{abstract}
For a family $\mathcal{F}$ of graphs, let $\ex(n,\mathcal{F})$ denote the maximum number of edges in an $n$-vertex graph which contains none of the members of $\mathcal{F}$ as a subgraph.
A longstanding problem in extremal graph theory asks to determine the function $\ex(n,\{C_3,C_4\})$.
Here we give a new construction for dense graphs of girth at least five with arbitrary number of vertices,
providing the first improvement on the lower bound of $\ex(n,\{C_3,C_4\})$ since 1976.
As a corollary, this yields a negative answer to a problem in Chung-Graham \cite{CG}.
%as a strong version of a well-known conjecture of Erd\H{o}s.
\end{abstract}

\section{Introduction}
For a given family $\mathcal{F}$ of graphs,
throughout this note we denote $\ex(n,\mathcal{F})$ to be the maximum number of edges in an $n$-vertex graph which does not contain any member in $\mathcal{F}$ as its subgraph.
This number - often referred as the {\it extremal number} or {\it Tur\'an number} of $\mathcal{F}$ - is the main subject in the field of extremal graph theory (see \cite{FS-survey}).
One of the central, extremely challenging problems in this field asks for the determination of the extremal number $\ex(n,\{C_3,C_4\})$ of the family consisting of the triangle $C_3$ and the 4-cycle $C_4$,
whose study can be dated back to a paper of Erd\H{o}s \cite{E38} in 1938.

A relevant object is the {\it Zarankiewicz number} $z(n,C_4)$ of the 4-cycle,
that is, the maximum number of edges in an $n$-vertex bipartite graph without containing a 4-cycle.
It is well-known that $z(n,C_4)=(\frac{n}{2})^{3/2}+o(n^{3/2})$ (see \cite{DHS,FS-survey}).
To be more precise, there exists a constant $c>0$ such that for any positive integer $n$,
\begin{equation}\label{equ:Z}
\left(\frac{n}{2}\right)^{3/2}-cn^{4/3}\leq z(n,C_4)\leq \frac{n}{4}\left(\sqrt{2n-3}+1\right)\leq \left(\frac{n}{2}\right)^{3/2}+\frac{1}{4}n,
\end{equation}
where the lower bound follows from \cite{F96} and the upper bound can be found in \cite{KSV} (see its Proposition 1.4).\footnote{Using the result of \cite{BHP} on the distribution of primes, the proof in \cite{F96} can yield a slightly better general lower bound that $z(n,C_4)\geq \left(\frac{n}{2}\right)^{3/2}-cn^{1.2625}$ for some $c>0$ and any positive integers $n$.}
As a bipartite graph cannot contain a triangle,
evidently one can relate these two aforementioned numbers with the following inequality
\begin{equation}\label{equ:ex>=z}
\ex(n,\{C_3,C_4\})\geq z(n,C_4).
\end{equation}
A famous old conjecture of Erd\H{o}s \cite{E38,E75} (restated in Erd\H{o}s-Simonovits \cite{ES}) asserts that this lower bound on $\ex(n,\{C_3,C_4\})$ essentially is optimal.

\begin{conjecture}[Erd\H{o}s \cite{E38,E75}, Erd\H{o}s-Simonovits \cite{ES}]\label{Conj:Erd}
It holds that
\begin{equation*}\label{equ:conj}
\lim\limits_{n\to +\infty} \frac{\ex(n,\{C_3,C_4\})}{z(n,C_4)}=1.
\end{equation*}
\end{conjecture}
In view of \eqref{equ:Z} and \eqref{equ:ex>=z}, this conjecture is equivalent to the upper bound $\ex(n,\{C_3,C_4\})\leq (\frac{n}{2})^{3/2}+o(n^{3/2})$.
It is still widely open.
The best known upper bound on $\ex(n,\{C_3,C_4\})$
remains the following trivial bound that
$$\ex(n,\{C_3,C_4\})\leq \ex(n,C_4)=\frac12 n^{3/2}+O(n).\footnote{Throughtout this note, for a function $f(n)$ we write $f(n)=O(n)$ if there exists some absolute constant $C>0$ such that $f(n)\leq Cn$ for any positive integers $n$.
Here, $\ex(n,C_4)=\frac12 n^{3/2}+O(n)$ is a well-known result of K\H{o}v\'ari-S\'os-Tur\'an \cite{KST} and Reiman \cite{R58}.}$$
On the other hand,
Parsons \cite{P76} gave a construction in 1976, showing that for integers $n=\binom{q}{2}$ where $q=1$ mod $4$ is a prime,
the inequality \eqref{equ:ex>=z} can be improved to
\begin{equation*}\label{equ:parson}
%\ex(n,\{C_3,C_4\})\geq z(n,C_4)+\frac{1}{8}n.
\ex(n,\{C_3,C_4\})\geq \left(\frac{n}{2}\right)^{3/2}+\frac{3}{8}n\geq z(n,C_4)+\frac{1}{8}n.
\end{equation*}
To the best of our knowledge, no progress has been made since then.
A stronger version of Conjecture~\ref{Conj:Erd} was stated as a problem in the book of Chung and Graham \cite{CG} (see p.41 in Section 3.4).
\begin{prob}[Chung-Graham \cite{CG}]\label{pro:strong-conj}
Is it true that $\ex(n,\{C_3,C_4\})=\left(\frac{n}{2}\right)^{3/2}+O(n)$?
\end{prob}
Let us also mention that Allen, Keevash, Sudakov and Verstra\"ete \cite{AKSV} (see Conjecture 1.7 therein) made an opposite conjecture that
$\liminf\limits_{n\to +\infty}\frac{\ex(n,\{C_3,C_4\})}{z(n,C_4)}>1.$

A general conjecture of Erd\H{o}s and Simonovits \cite{ES} concerning extremal numbers of families containing bipartite graphs and odd cycles is as follows:
Let $\mathcal{C}$ denote the family of all odd cycles, and let $\mathcal{C}_k$ denote the family of all odd cycles of length at most $k$.
Then for any finite family $\mathcal{F}$ containing a bipartite graph, there exists an odd integer $k$ such that
$$\lim\limits_{n\to +\infty} \frac{\ex(n,\mathcal{F}\cup \mathcal{C}_k)}{\ex(n,\mathcal{F}\cup \mathcal{C})}=1.$$
Erd\H{o}s and Simonovits \cite{ES} confirmed this for $\mathcal{F}=\{C_4\}$ by showing that $\ex(n,\{C_4,C_5\})=\left(\frac{n}{2}\right)^{3/2}+O(n)$
and thus $\lim\limits_{n\to +\infty}\frac{\ex(n,\{C_4,C_5\})}{z(n,C_4)}=1$.
This was strengthened by Keevash, Sudakov and Verstra\"ete \cite{KSV}, where their main result implies that for all integers $k\geq 2$,
\begin{equation}\label{equ:C4+C2k+1}
\ex(n,\{C_4,C_{2k+1}\})=\left(\frac{n}{2}\right)^{3/2}+O(n).
\end{equation}

The main result of this paper is the following theorem,
which improves the error term in Parsons' lower bound on $\ex(n,\{C_3,C_4\})$ from $\Omega(n)$ to $\Omega(n^{1.25})$.

\begin{thm}\label{thm:main}
There exists an absolute constant $c>0$ such that for every integer $n\geq 7$,
$$\ex(n,\{C_3,C_4\})\geq  z(n,C_4)+c\cdot n^{1.25}.$$
\end{thm}

We would like to emphasize that this result works for {\it every} integer $n\geq 7$,
while the construction of Parsons is applicable only for a special form of integers $n$.
Let us also note that when $n=6$, both numbers $\ex(n,\{C_3,C_4\})$ and $z(n,C_4)$ are equal to $6$.

Using the above bound, one can immediately derive the following corollary that 
the difference between $\ex(n,\{C_3,C_4\})$ and $\left(\frac{n}{2}\right)^{3/2}$ can be a superlinear term in $n$,
thus showing that Problem~\ref{pro:strong-conj} does not hold in general.

\begin{corollary}\label{coro}
For integers $n=2(q^2+q+1)$ where $q$ is a prime power,
$$\ex(n,\{C_3,C_4\})=\left(\frac{n}{2}\right)^{3/2}+\Omega(n^{1.25}).$$
In particular, this provides a negative answer to Problem~\ref{pro:strong-conj}.
\end{corollary}

We also see from \eqref{equ:C4+C2k+1} that the behavior of $\ex(n,\{C_3,C_4\})$ 
is different from $\ex(n,\{C_4,C_{2k+1}\})$ for any integer $k\geq 2$, as they vary in their second order terms.
We remark (see the last paragraph of Section~2) that using results from number theory on the distribution of primes,
the conclusion of Corollary~\ref{coro} in fact holds for almost all integers $n$.

The rest of the paper is organized as follows.
In Section~2, we present the proofs of Theorem~\ref{thm:main} and Corollary~\ref{coro}.
In the final section, we conclude with a remark by explaining that similar constructions as in Theorem~\ref{thm:main} are unlikely to give better bounds.

\section{The proofs of Theorem~\ref{thm:main} and Corollary~\ref{coro}}

In this section, we first prove Theorem~\ref{thm:main} and then use it to infer Corollary~\ref{coro}.

\begin{proof}[\bf Proof of Theorem \ref{thm:main}]
We begin with a warm-up by showing that for any integer $n\geq 7$,
\begin{equation}\label{equ:+1}
\ex(n,\{C_3,C_4\})\geq z(n,C_4)+1.
\end{equation}
By \eqref{equ:ex>=z} we have $\ex(n,\{C_3,C_4\})\geq z(n,C_4)$.
Assume for a contradiction that $\ex(n,\{C_3,C_4\})= z(n,C_4)$ for some $n\geq 7$.
Then there exists an $n$-vertex bipartite $C_4$-free\footnote{Given a family $\mathcal{F}$ of graphs, throughout the rest we say a graph $H$ is {\it $\mathcal{F}$-free} if $H$ does not contain any member in $\mathcal{F}$ as a subgraph. If $\mathcal{F}$ consists of a single graph $F$, then simply we say {\it $F$-free} instead of $\{F\}$-free.} graph $G$ with $e(G)=\ex(n,\{C_3,C_4\})$.
Let $(X,Y)$ be the bipartition of $G$ with $|X|\geq |Y|$.
We claim that
\begin{itemize}
	\item[(P1).] Any two vertices in the same part ($X$ or $Y$) have a unique common neighbor, and
	\item[(P2).] The maximum degree $\Delta(G)$ is at most three.
\end{itemize}
If $u,v\in X$ do not have a common neighbor, then the graph $G+\{uv\}$ will be $\{C_3,C_4\}$-free,
a contradiction to that $e(G)=\ex(n,\{C_3,C_4\})$.
So (P1) follows.
To see (P2), suppose that there is a vertex $u\in X$ with at least four neighbors $a,b,c,d\in Y$.
Let $G'$ be the graph obtained from $G$ by deleting the edges $ub,uc$ and adding new edges $ab,bc,cd$.
%Note that $\{u,a,b,c,d\}$ induces a cycle of length five in $G'$. and one can easily verify that $G'$ is also $\{C_3,C_4\}$-free with $e(G')>e(G)$.
If there is a  3-cycle or 4-cycle in $G'$ containing one edge in $\{ab,bc,cd\}$, then this cycle contains a vertex $u'\in X$ which has two neighbors in $\{a,b,c,d\}$. It implies that $u$ and $u'$ have two common neighbors in $G$, which is impossible.   Thus $G'$ is   $\{C_3,C_4\}$-free with $e(G')>e(G)$.
This contradiction proves the above claim.
Let $\sigma_X$ denote the number of paths of length two with both end-points in $X$.
By (P1) we have $\sigma_X=\binom{|X|}{2}$, while (P2) implies that $\sigma_X=\sum_{y\in Y}\binom{d(y)}{2}\leq 3|Y|$.
As $|X|\geq  |Y|$, we get $\binom{|X|}{2}\leq 3|Y|\leq 3|X|$, showing that $|X|\leq 7$.
So $n\leq 14$.
The precise values of $\ex(n, \{C_3,C_4\})$ are determined in \cite{GKL} for all integers $n\leq 24$ (see Theorem~3.1).
In particular, when $n\in \{7,8,9,10,13\}$,
the extremal $\{C_3,C_4\}$-free graph on $n$ vertices is unique and non-bipartite (see Figures 1 and 2 in \cite{GKL}).
So it only remains to consider $n\in \{11,12,14\}$.
For $n=11$, we have $|Y|\leq 5$ and $\ex(11,\{C_3,C_4\})=16$ from \cite{GKL},
so there must be a vertex in $Y$ of degree at least four, a contradiction to (P2).
For $n=12$, we have $\ex(12,\{C_3,C_4\})=18$ from \cite{GKL}, implying that all vertices have degree three and $|X|=|Y|=6$.
However, it leads to a contradiction as we should have $\binom{|X|}{2}=\sigma_X=3|Y|$ in this situation.
Lastly, for $n=14$, we have $\ex(14,\{C_3,C_4\})=23$ from \cite{GKL} and
thus the extremal graph $G$ must contain a vertex of degree at least four, a contradiction to (P2). This proves \eqref{equ:+1}.

To complete the proof of Theorem~\ref{thm:main}, it suffices to show that for sufficiently large $n$,
\begin{equation}\label{equ:Omega}
\ex(n,\{C_3,C_4\})=z(n,C_4)+\Omega(n^{1.25}).
\end{equation}
Let $\varepsilon$ be a sufficiently small (but fixed) positive real and let $n$ be any integer which is sufficiently larger than $1/\epsilon$.
Let $G$ be any extremal graph of $z(n,C_4)$, i.e., an $n$-vertex bipartite $C_4$-free graph with $z(n,C_4)$ edges.
Let $(X,Y)$ be the bipartition of $G$.
In what follows, based on $G$ we will construct a (non-bipartite) $\{C_3,C_4\}$-free graph on the same vertex set $V(G)$ and with $\Omega(n^{1.25})$ more edges.

We claim that there exists a vertex $u$ in $G$ with $d(u)\leq (1+\varepsilon)\sqrt{n/2}$ and $|\cup_{x\in N(u)} N(x)|\geq (1-\varepsilon)n/2$.
Note that as $G$ is $C_4$-free, for any vertex $u$, the neighborhoods $N(x)$'s for all vertices $x\in N(u)$ are pairwise disjoint.
To prove this claim, we will proceed to show that
\begin{itemize}
\item [(A).] there are less than $n/2$ vertices with degree at least $(1+\varepsilon)\sqrt{ n/2}$, and
\item [(B).] there are less than $n/2$ vertices $u$ with $|\cup_{x\in N(u)} N(x)|\le (1-\varepsilon) n/2$.
\end{itemize}
First let us see that $(X,Y)$ is almost balanced.
By \eqref{equ:Z}, we have
$$\left(\frac{n}{2}\right)^{3/2}-cn^{4/3}\leq z(n,C_4)=e(G)\leq (|X||Y|)^{3/4}+\max\{|X|,|Y|\}\leq (|X|(n-|X|))^{3/4}+n,$$
where the second last inequality follows by Proposition 3.9 in \cite{KSV}.
Solving the above inequality for $|X|$, it gives that $(1-\frac{\varepsilon}{2})\frac{n}{2}\leq |X|, |Y|\leq (1+\frac{\varepsilon}{2})\frac{n}{2}$.
Let $\sigma$ be the number of paths of length two in $G$.
As $G$ is bipartite and $C_4$-free,
each pair of vertices from the same part is contained in at most one path of length two.
So $\sigma\leq \binom{|X|}{2}+\binom{|Y|}{2}\leq (1+\frac{\varepsilon^2}{4})\frac{n^2}{4}$.
Suppose for a contradiction to (A) that there are $n/2$ vertices with degree at least $(1+\varepsilon)\sqrt{n/2}$ in $G$.
By Jensen's inequality, we have
\begin{align*}
\sigma &=\sum_{u\in V(G)}\binom{d(u)}{2}
\geq \frac n2 \binom{(1+\varepsilon)\sqrt{n/2}}{2}+ \frac n2 \binom{\left(2e(G)-\frac n2(1+\varepsilon\right)\sqrt {n/2})/(n/2)}{2}\\
&\ge \frac n2 \binom{(1+\varepsilon)\sqrt{n/2}}{2}+ \frac n2 \binom{(1-\varepsilon)\sqrt{n/2}-O(n^{1/3})}{2}=(1+\varepsilon^2)\frac{n^2}{4}-o(n^2).
\end{align*}
This is a contradiction to the above upper bound on $\sigma$, thus proving (A).
To see (B), suppose on the contrary that there are $n/2$ vertices $u$ with $|\cup_{x\in N(u)} N(x)|\le (1-\varepsilon) n/2$.
Each of these vertices is contained in at most $(1-\varepsilon) n/2$ paths of length two as an end-point,
while any other vertex is contained in at most $\max\{|X|,|Y|\}\leq (1+\frac{\varepsilon}{2})\frac{n}{2}$ paths of length two as an end-point.
Totally we have
$2\sigma \leq \frac{n}{2}\cdot (1-\varepsilon)\frac{n}{2}+ \frac{n}{2}\cdot(1+\frac{\varepsilon}{2})\frac{n}{2}=\left(2-\frac{\varepsilon}{2}\right)\frac{n^2}{4},$
implying that
\begin{align*}
\left(1-\frac{\varepsilon}{4}\right)\frac{n^2}{4}\geq \sigma =\sum_{u\in V(G)}\binom{d(u)}{2}
\geq n \binom{ 2e(G)/n}{2}\geq n \binom{\sqrt{n/2}-O(n^{1/3})}{2}=\frac{n^2}{4}-o(n^2),
\end{align*}
a contradiction. This completes the proof of the claim.

Let $u\in X$ be the vertex as claimed.
Let $N(u)=\{u_1,u_2,...,u_t\}$ for some $t\leq (1+\varepsilon)\sqrt{n/2}$.
For each $1\leq i\leq t$, let $N_i=N(u_i)$ and $E_i=\{u_ix| x\in N_i\}$.
As pointed out, these $N_i$'s are pairwise disjoint and thus $\sum_{1\leq i\leq t} |N_i|\geq (1-\varepsilon)n/2.$
Let $G_i$ be an extremal $\{C_3,C_4\}$-free graph on the vertex set $N_i$.
By \eqref{equ:Z}, there exists some $c>0$ such that
\begin{equation}\label{equ:e(Gi)}
e(G_i)\geq \left(|N_i|/2\right)^{3/2}-c\cdot\left(|N_i|/2\right)^{4/3}.
\end{equation}

Let $H$ be obtained from $G$ by deleting all edges in $E_i$ and adding the graph $G_i$ into $N_i$ for every $1\leq i\leq t$.
We claim that $H$ is $\{C_3, C_4\}$-free.
Suppose that $H$ contains a triangle say $abc$.
Then at least one edge (say $ab$) must appear in some $G_i$ (note that $V(G_i)=N_i\subseteq X$).
As $G_i$ is $C_3$-free, we must have $c\in Y$ and thus $ac,bc\in E(G)\cap E(H)$.
But this is a contradiction as the unique common neighbor of $a,b$ in $Y$ has been destroyed by deleting the edges of $E_i$.
Now suppose $H$ has a $C_4$ say $abcd$. We may assume that $ab\in E(G_i)$.
Since $G_i$ is $C_4$-free and $H[Y]$ is an independent set,
we may assume that $c\in X$ and $d\in Y$.
As $N_j$'s are pairwise disjoint, it is clear that $c\in N_i$.
Then we get $ad,cd\in E(G)\cap E(H)$, which is a contradiction by the same reason.
So indeed, $H$ is an $n$-vertex $\{C_3, C_4\}$-free graph.
We can estimate the number of edges in $H$ as follows:
\begin{align*}
e(H)&=e(G)+\sum_{1\leq i\leq t} (e(G_i)- |E_i|)\geq z(n,C_4)+ \sum_{1\leq i\leq t} \left(\left(|N_i|/2\right)^{3/2}-c\cdot\left(|N_i|/2\right)^{4/3}-|N_i|\right)\\
&\geq z(n,C_4)+\sum_{1\leq i\leq t} (1-o(1))\left(\frac{|N_i|}{2}\right)^{3/2}
\ge z(n,C_4)+ (1-o(1))\cdot t\left(\frac{\sum_{1\leq i\leq t}|N_i| /t}{2}\right)^{3/2}\\
&%=z(n,C_4)+ (1-o(1)) 2^{-1.5} d(u)  ^{-0.5} \left( \sum_{u_i\in N(u)}|N_i| \right)^{1.5}
=z(n,C_4)+ \Omega\left(\frac{(\sum_{1\leq i\leq t}|N_i|)^{3/2}}{\sqrt{t}}\right)\geq z(n,C_4)+\Omega(n^{1.25}),
\end{align*}
where the first inequality follows by \eqref{equ:e(Gi)}, and the last inequality uses the facts that $t\leq (1+\varepsilon)\sqrt{n/2}$ and $\sum_{1\leq i\leq t} |N_i|\geq (1-\varepsilon)n/2.$
This proves \eqref{equ:Omega} and thus completes the proof of Theorem~\ref{thm:main}.
\end{proof}

%One can derive Corollary~\ref{coro} promptly.

\begin{proof}[\bf Proof of Corollary~\ref{coro}]
Let $q$ be a prime power and let $n=2(q^2+q+1)$.
%In this case, it is well-known that a {\it generalized 3-gon} of order $q$ exists
In this case, it is well-known that a finite projective plane of order $q$ exists
and thus $z(n,C_4)=\frac12(q+1)n\geq \left(\frac{n}{2}\right)^{3/2}$ (see e.g. Theorem~1.2 in \cite{KSV}).
Therefore, by Theorem~\ref{thm:main}, we have
$\ex(n,\{C_3,C_4\})=z(n,C_4)+\Omega(n^{1.25})=\left(\frac{n}{2}\right)^{3/2}+\Omega(n^{1.25}).$
\end{proof}

In the rest of this section, we present a generalization of Corollary~\ref{coro}.
Let $p_n$ denote the $n$'{th} prime.
It is well-known in number theory (see \cite{HB}) that
there exists some constant $\delta\in (0,1)$ such that for all reals $x>0$,
\begin{equation}\label{equ:pn}
\sum_{p_n\leq x, ~ p_{n+1}-p_n\geq \sqrt{p_n}} (p_{n+1}-p_n)=O(x^{1-\delta}).
\end{equation}
This can imply that for sufficiently small constant $\epsilon>0$ and for almost all integers $n$, there exists a prime in $[n-\epsilon\sqrt{n},n]$.
By the proof of \cite{F96}, one can then derive that $z(n,C_4)\geq \left(\frac{n}{2}\right)^{3/2}-O(\epsilon)\cdot n^{1.25}$ for almost all integers $n$.
Therefore, together with Theorem~\ref{thm:main}, this shows that
$\ex(n,\{C_3,C_4\})=z(n,C_4)+\Omega(n^{1.25})=\left(\frac{n}{2}\right)^{3/2}+\Omega(n^{1.25})$ holds for almost all integers $n$.

\section{A remark}
In the above construction, we take an extremal graph of $z(n,C_4)$, choose vertex-disjoint subsets of size roughly $q=\sqrt{n/2}$,
and then for each of these subsets $A$, add $\Omega(q^{3/2})$ edges into $A$ and delete fewer edges incident with $A$ to make a $\{C_3, C_4\}$-free graph with more than $z(n,C_4)$ edges.
One may ask whether one can take larger subsets (say of size $n^{1/2+\epsilon}$ for any $\epsilon>0$) and add/delete edges using similar operations to get a denser $\{C_3, C_4\}$-free graph.
We illustrate in the following example that it is unlikely to give better constructions.

For the purpose of our presentation,
let $G$ be an extremal graph of $z(n,C_4)$ with bipartite $(X,Y)$ and let $n=2q^2$ be an integer with $q\in \mathbb{R}^+$ such that
\begin{equation}\label{equ:MYconj}
\mbox{except $O(q)$ vertices, every vertex in $G$ has degree at least $q-o(q)$.}
\end{equation}
Let $\delta\in (0,1)$ be any real. Consider any set $A\subseteq X$ of size $q^{1+\delta}$.

We will show that it is impossible to construct a $\{C_3,C_4\}$-free graph $H$ obtained from $G$ by adding $c |A|^{3/2}$ edges into $A$ and deleting any subset $E^*$ of edges such that $e(H)>e(G)$.
Suppose for a contradiction that such $H$ does exist.
Since $e(H)>e(G)$, we have $|E^*|<c |A|^{3/2}=o(q^{2+\delta})$.
It is easy to see that the size of $|X|$ or $|Y|$ is $(1+o(1))q^2$.
By \eqref{equ:MYconj},
the number of edges between $A$ and $Y$ in $G$ is at least $(1-o(1))q^{2+\delta}$.
Then the induced bipartite subgraph $H[A,Y]$ of $H$ with parts $A$ and $Y$ has at least $(1-o(1))q^{2+\delta}-|E^*|\geq (1-o(1))q^{2+\delta}$ edges.
Let $\sigma_1$ be the number of paths of lengths two in $H[A,Y]$ with both ends in $A$. Then
\begin{align*}
\sigma_1=\sum_{v\in Y}\binom{d_{H[A,Y]}(v)}{2}\geq |Y|\binom{e(H[A,Y])/|Y|}{2}=(1-o(1))q^{2+2\delta}/2=(1-o(1))\binom{|A|}{2}.
\end{align*}
As $H[A]$ has at least $c|A|^{3/2}$ edges,
the number $\sigma_2$ of paths of length two in $H[A]$ is
\begin{align*}
\sigma_2=\sum_{v\in A}\binom{d_{H[A]}(v)}{2}\geq |A|\binom{2e(H[A])/|A|}{2}\geq c^2\binom{|A|}{2},
\end{align*}
Therefore, in total $H$ contains $\sigma_1+\sigma_2\geq (1+c^2-o(1))\binom{|A|}{2}>\binom{|A|}{2}$ paths of length two with both ends in $A$,
which leads to a copy of $C_4$ in $H$, a contradiction.

We point out that when $n=2(t^2+t+1)$ for any prime power $t$, the extremal graph of $z(n, C_4)$ is regular and thus satisfies \eqref{equ:MYconj},
so (at least) for these infinitely many integers $n$, our construction cannot be improved using the above operations.
In fact, to make the above arguments work, all we need here is the property that $e(G[A,Y])\geq (1-o(1))q^{2+\delta}$ (which follows by \eqref{equ:MYconj}),
and one can show that for general $n$, almost all subsets $A$ of size $q^{1+\delta}$ in one part of the extremal graph of $z(n,C_4)$ satisfy this property.

As a side note, Keevash et al. proved in \cite{KSV} (see Theorem~5.1) that any (nearly) extremal graph of $z(n,C_4)$ satisfies the pseudorandomness property.
We would like to conjecture that any extremal graph of $z(n,C_4)$ satisfies \eqref{equ:MYconj}.

\bigskip

\noindent {\bf Acknowledgements.}
We would like to thank Lilu Zhao for his helpful instruction in number theory.
We are also grateful to Chunqiu Fang for the stimulating discussion and Guorong Gao for bringing the result on \eqref{equ:pn} to our attention.

\bibliographystyle{unsrt}

\end{document}